\documentclass[reqno]{amsart}
\usepackage{amsmath}
\usepackage{amssymb}
\usepackage{stmaryrd}
\title{
A Simple Regularization of Graphs
}
\author{Yoshiyasu Ishigami}
\address{Department of Information and Communication Engineering,
The University of Electro-Communications, Chofu, Tokyo 182-8585, Japan.}
\date{}
\def\picture #1 by #2 (#3){
                \vbox to #2{
                        \hrule width #1  height 0pt depth
0pt
                        \vfill
                        \special{picture #3}}}
\def\scaledpicture #1 by #2(#3 scaled #4){{
                \dimen0=#1 \dimen1=#2
                \divide\dimen0 by 1000 \multiply\dimen0 by
#4
                \divide\dimen1 by 1000 \multiply\dimen1 by
#4
                \picture\dimen0 by \dimen1 (#3 scaled #4)}}
\newtheorem{adf}{Definition}[section]

\newtheorem{thm}{Theorem}[section]
\newtheorem{cor}[thm]{Corollary}
\newtheorem{lemma}[thm]{Lemma}

\newtheorem{fact}[thm]{Fact}

\newtheorem{remark0}[thm]{Remark}

\newtheorem{algorithm0}[thm]{Algorithm}
\newtheorem{asett}[thm]{Setup}

\newenvironment{df}{
\begin{adf}\begin{sl}}
{\end{sl}\thqed\end{adf}}
\newenvironment{sett}{
\begin{asett}\begin{sl}}{\end{sl}\thqed\end{asett}}

\newcommand{\thqed}{\hfill\fbox{}\\ }
\newcommand{\Proof }{{\bf Proof} : }

\newcommand{\factqed}{\hfill\bigskip\fbox{}\\}

\newcommand{\lemmaqed}{\hfill\bigskip\fbox{}\\}

\def\brkt#1{\left({#1}\right)}
\def\dbrkt#1{\llbracket{#1}\rrbracket}

\def\ang#1{\langle{#1}\rangle}

\newcommand{\naturalset}{{\mathbb N}}

\newcommand{\Prob}{{\mathbb P}}
\newcommand{\Ex}{{\mathbb E}}

\newtheorem{procedure0}{Procedure}

\newcount\timehh\newcount\timemm\timehh=\time\divide\timehh by 60
\timemm=\time\count255=\timehh\multiply\count255 by-60
\advance\timemm by \count255
\ifnum\timehh=12\def\apm{pm}\else
\ifnum\timehh>12\def\apm{pm}\advance\timehh by-12\else
\def\apm{am}\fi\fi
\def\timestamp{\number\timehh\,:\,\ifnum\timemm<10
0\fi\number\timemm\,\apm}
\begin{document}
\maketitle
\begin{abstract}
The well-known regularity lemma of E.~Szemer\'edi for graphs 
(i.e. $2$-uniform hypergraphs) claims that for any graph 
there exists a vertex partition with the property of quasi-randomness.
We give a simple construction of such a partition.
It is done just by 
taking a constant-bounded number of random vertex samplings only one time (thus, iteration-free).
Since it is independent from the definition of quasi-randomness, 
it can be generalized very naturally to hypergraph regularization. 
In this expository note, we show only a graph case of the paper \cite{I06} on 
hypergraphs, 
but may help the reader to access \cite{I06}.
\end{abstract}
\section{Introduction}
The well-known regularity lemma of Szemer\'edi \cite{Sz} (also called the uniformity lemma) 
was discovered in the course of 
obtaining the so-called Szemer\'edi's theorem on arithmetic progressions \cite{Sz75} as 
an affirmative answer of a conjecture by Erd\H{o}s and Tur\'an.
 It has been known that this graph-theoretic lemma has a plenty of applications 
 in many topics of mathematics and 
theoretical computer sciences. 
\par
The regularity lemma claims that for any ordinary graph (i.e. any 2-uniform hypergraph)
there exists a vertex partition with the property of quasi-randomness. 
Our purpose of this note is to give a simple construction of such a partition. 
It has several advantages over previously-known methods.
It is the case of $k=2$ (i.e. the case of 2-uniform hypergraphs) in \cite{I06} which deals with general $k$. 
Although this expository note is not necessary to read \cite{I06}, skimming it 
may help the reader understand the main idea of \cite{I06}.\par
Remark that our construction had not been known even for the simplest case $k=2$ before \cite{I06}.
Although the idea of partitioning the vertex-set randomly has been previously known 
(\cite{GGR, AFNS}), such a construction was done always by a constant number of sample random vertices, 
which thus needs an iteration. 
The key difference is that ours is constant-bounded but the number of random samplings is chosen also randomly.
The proof can be naturally deduced once the claim is given.
\par
Recall how the standard proof by Szemer\'edi constructs the desired vertex partition with 
quasi-randomness. 
Roughly speaking, the partition was constructed by iterated 
applications of the dichotomy between energy-increment and structure. 
That is, initially take an arbitrary vertex partiton (with a constant number of vertex sets). 
It can be shown that 
\\
(1) this partition satisfies the required quasi-random property 
or that
\\
 (2) there must exist another vertex 
partition finer than this partition such that 
\par
(2.1) the number of vertex sets increases 
but is still bounded by a constant and further that 
\par
(2.2)  a value called \lq energy\rq\ 
(or \lq index\rq\ ) of the finer partition is significantly larger than the 
\lq energy\rq\ of the coarser partition. \\
They replace the coarser partition by the finer one and repeat this process.
Since the energy is always less than one from its definition, 
the repeating process must stop in at most constant time.
(Note that however the exact time when it stops 
depends on the structure of the given graph.)
The vertex partition  which the final stage outputs satsifies (1) and thus is the desired partition.
\par
On the other hand, our construction goes as follows.
\\
(0')
Take a large constant $\tilde{n}$ which depends on $\epsilon$ (parameter on how much quasi-random it 
should be) but 
is independent from (the number of vertices of) the given graph.
Further take a length-$\tilde{n}$
 integer sequence $0=m_0\ll m_1\ll \cdots \ll 
m_{\tilde{n}-1},$ also independent from the given graph.
\\
(1')
Choose an integer $0\le n<\tilde{n}$ uniformly at random and 
further choose $m_n$ vertices 
uniformly at random from the given graph.
\\
(2')
Each vertex of the given graph is labeled 
by the adjacency between the vertex and the randomly-chosen vertices.
\\
The resulting partition certainly consists of a constant number of vertex sets (i.e. $2^{m_n}\le 
2^{m_{\tilde{n}-1}}$) 
and would be the desired partition with high probability.
\par
Previous constructions including the usual one by Szemer\'edi 
consist of iterated procedures, 
while our construction consists of only one procedure.
Furthermore ours is independent from the definition of quasi-randomness, 
while previous constructions depend on it.
Several definitions of quasi-randomness have been known. 
For the case of ordinary graphs, all of them are known to be equivalent.
However, it has been noted that unlike the situation for graphs, there are
several ways one might define regularity for hypergraphs
(Tao-Vu \cite[pp.455]{TV},R\"{o}dl-Skokan \cite[pp.1]{RSk04}).
Because our construction is independent from the definition of 
quasi-randomness, it can be naturally generalized from ordinary graphs 
to hypergraphs.
For the extension to hypergraphs, see \cite{I06}.
\bigskip
\par
The purpose of this note is to present a new construction of 
the vertex partition (for the case of ordinary graphs) 
and to show that it certainly satisfies quasi-randomness.
For the case of ordinary graphs, there are several definitions of 
quasi-randomness but all of them are known to be equivalent (\cite{CGW}
).
As our definition of quasi-randomness, we will choose the number of 
induced subgraphs. 
In this paragraph, I will explain why we will use this definition, 
though of course it is not serious at least when considering only 
the case of ordinary graphs (since they are equivalent).
The usual regularity lemma firstly defines quasi-randomness by 
a condition on the number of edges between two sets of vertices.
Secondly Szemer\'edi proved the existence of a vertex partition with 
this quasi-randomness. Thirdly and finally, 
the quasi-randomness on counting induced subgraphs (i.e. our definition) 
can be derived from his quasi-randomness (on edges between two subsets).
The third step is called to be 
the {\em counting lemma} and is easy to show, so the second step only is 
the core of the matter.
All of the three hypergraph-theoretic 
proofs of Szemer\'edi's arithmetic-progression theorem 
by R\"{o}dl et al. \cite{RSk04,NRS}, Gowers \cite{G} and Tao \cite{Tao06} 
can be considered as generalizations of the above three steps.
But unlike the case of graphs, the third step (counting lemma) 
was hard to show for hypergraphs. All of the three proofs are different 
partly because all of them employed different 
definitions of quasi-randomness, on which their regularizations depend.
On the other hand, we will not follow the above three steps. 
We will define a probabilistic construction for partitioning 
the vertices, 
which will be proven to satisfy the condition of our quasi-randomness 
on counting induced subgraphs. This strategy can 
 be very naturally generalized from graphs to hypergraphs in \cite{I06}.
I believe that this framework of hypergraph regularity 
lemma is convenient for 
a wide range of applications on hypergraphs. 
In fact, applications of our method are seen 
in \cite{I06m,I06lr,I07}. 
\bigskip
\par
One of the new major technical ingredients in our proof 
comes 
from the use of \lq linearity of expectation.\rq\ All of the
previous proofs use the dichotomy (or energy-increment)
explicitly. (See \cite[\S6]{G},
 \cite[\S1]{Tao06a}.) Namely, when proving 
the existence of a vertex partition, they define an \lq
energy\rq\ (or index)
by the {\em maximum} (or supremum) of some (energy) function. (For
example, see \cite[eq.~(8)]{Tao06}.) It corresponds to
(\ref{061220}) in this paper. They consider the {\em maximum}
value of this energy over all subdivisions in each step. If the
energy significantly increases by {\em some} subdivision, they
take the {\em worst} subdivision as the base partition of the next
step. They then repeat this process. Since the energy is bounded,
this operation must stop at some step, in which case
 there is no quite bad subdivision, and thus, most cells should be
quasi-random
(dichotomy).
\par
On the other hand, we (implicitly) take an {\em average}
subdivision instead of the worst one. The definition of our
regularization determines the probability space of partitions
(subdivisions). We also randomly decide on the number of vertex
samples to choose.

With these ideas, we can {\em hide} the troublesome dichotomy
iterations inside linear equations of expectations (\ref{lj20}).
\par
\bigskip
\par
We have two reasons why we will deal with multi-colored graphs
instead of ordinary graphs, even though almost all previous
researchers dealt with the usual graphs.
 First, our proof of the regularity lemma will be natural.
Second, we can naturally combine subgraph (black\&invisible) and
induced-subgraph (black\&white) problems when we apply our result,
while the two have usually been discussed separately.

\allowdisplaybreaks[4]
\section{Statement of the Theorem}
In this paper, $\Prob$ and $\Ex$ will denote probability and
expectation, respectively. We denote conditional probability and
expectation by $\Prob[\cdots|\cdots]$ and $\Ex[\cdots|\cdots].$
\begin{sett}\label{f1231}
Throughout this paper, we fix a positive integer $r$ and an \lq
index\rq\ set $\mathfrak{r}$ with $|\mathfrak{r}|=r.$ Also we fix
a probability space $({\bf \Omega}_i,{\mathcal B}_i,\Prob)$ for
each $i\in \mathfrak{r}$. We assume that ${\bf \Omega}_i$ is
finite and that ${\mathcal B}_i=2^{{\bf \Omega}_i}$ (for the sake
of simplicity). Write ${\bf \Omega}:=({\bf \Omega}_i)_{i\in
\mathfrak{r}}$.
\end{sett}
In order to avoid using measure-theoretic jargon such as
measurability or Fubini's theorem, for the benefit of readers who
are interested only in applications to discrete mathematics, we
assume ${\bf \Omega}_i$ to be a (non-empty) finite set. However,
our arguments should be extendable to a general probability space.
For applications, ${\bf \Omega}_i$ usually will contain a huge
number of vertices. We will not use this assumption 
logically in our proof, but it will be important in our theorems that 
some parameters and functions depend on $r$ but independent from any $|{\bf \Omega}_i|.$
\par
In what follows, we will try to embed a small $r$-partite graph to another large $r$-partite 
graph, where the $r$ vertex sets of the large graph will be always $({\bf \Omega}_i)_{i\in\mathfrak{r}}.$
And the large graph and its vertices and edges will be denoted by bold fonts (ex. ${\bf G}, {\bf v}, 
{\bf v'}, 
{\bf e}, \cdots$).
\par
For an integer $a$, we write $[a]:=\{1,2,\cdots,a\},$ and
${\mathfrak{r}\choose [a]}:=\dot{\bigcup}_{i\in
[a]}{\mathfrak{r}\choose i} =\dot{\bigcup}_{i\in [a]}\{I\subset
\mathfrak{r}| |I|=i\}.$ Thus 
${\mathfrak{r}\choose [2]}
={\mathfrak{r}\choose 
1}\dot{\cup}{\mathfrak{r}\choose 2}=\mathfrak{r}\dot{\cup}{\mathfrak{r}\choose 2}.$
When $r$ disjoint sets $X_i, i\in {\mathfrak r},$
with indices from ${\mathfrak r}$ are called {\bf vertex sets}, we
write $X_J:=\{Y\subset \dot{\bigcup}_{i\in J}X_i\,| |Y\cap
X_j|=1\,, \forall j\in J\}$ whenever $J\subset {\mathfrak r}$.
Thus $|X_J|=\prod_{j\in J}|X_j|.$ That is, for $J=\{1,2\}$, $|X_J|=|X_1||X_2|.$
\begin{df}[Colored graphs]\label{070724}
\quad Suppose Setup \ref{f1231}. Given $b_1$ and $b_2,$ a {\bf $(b_1,b_2)
$-colored ($\mathfrak{r}$-partite) graph} $H$ is a triple
$((X_i)_{i\in\mathfrak{r}},({C}_I)_{I\in {\mathfrak{r}\choose
[2]}}, (\gamma_I)_{I\in {\mathfrak{r}\choose [2]}} )$ where:\\ (1)
each $X_i$ is a set called a \lq vertex set,\rq\ \\(2) ${C}_I$ is
a set with at most $b_{|I|}$ elements, and \\(3) $\gamma_I$ is a
map from $X_I$ to ${C}_I.$\\
We write $V(H)=\dot{\bigcup}_{i\in \mathfrak{r}}X_i
$ and ${\rm C}_I(H)={C}_I$ for $
I.$
Each element of $V(H)$ is called a {\bf vertex}. Each element
$e\in V_I(H)=X_I, I\in {\mathfrak{r}\choose [2]},
$ is called
 an {\bf (index-$I$) edge}. Thus, when $|I| = 1$, an
 index-$I$ edge is just a vertex of $H$. Each member in
 ${\rm C}_I(H)$ is a {\bf
(face-)color (of index $I$)}. Write $H(e)=\gamma_I(e)$ for each
$I.$ (So we will not need the notation $\gamma$ after this definition.)
\par
When $I=\{i,j\}\in {\mathfrak{r}\choose 2} $ (i.e. $i\not=j$) and $e=\{v_i,v_j\}\in V_I(H),$
we define the {\bf frame-color} and {\bf
total-color} of $e$ by ${H}(\partial
e):=({H}(v_i),H(v_j))$ and by
${H}(\ang{e})=H\ang{e}:=(
H(e); H(v_i),H(v_j)).$
For a vertex $v_i\in X_i$ (which is also an index-$i$ edge),
we define the total-color of $v$ by
$H(\ang{v})=H\ang{v}:=H(v).$
The frame-color of a vertex is always the empty $()$.
Write ${\rm TC}_I(H):=\{H\ang{e}|\,{e}\in X_I=V_I(H)\},$ ${\rm
TC}_s(H):=\bigcup_{I\in {\mathfrak{r}\choose s}}{\rm TC}_I(H),$
and ${\rm TC}(H):={\rm TC}_1(H)\dot{\cup}{\rm TC}_2(H),$ where ${\rm TC}$ means total-color.
\end{df}
As usual, we will call a $(b_1,b_2)$-colored graph just a {\bf colored graph} or a {\bf graph} 
when we do not need to mention values $b_1, b_2.$
\begin{df}[Complexes]
A {\bf (simplicial-)complex} is a (colored ${\mathfrak r}$-partite) graph
such that: \\
(1) for each $I\in {\mathfrak{r}\choose [2]}$ there exists
at most one index-$I$ color called \lq invisible\rq\ and that \\
(2) if 
(the color of) an edge $e$ is invisible then for any edge $e^*\supset
e,$ its color must also be invisible. \\
A color is {\bf visible} if and only if it is not invisible. We simply say that 
an edge is visible/invisible when its color is so.
\par
For a graph ${\bf G}$ on ${\bf \Omega}$, let ${\mathcal S}_{h,{\bf
G}}$ be the set of complexes $S$ such that:\\ 
(1) each of $r$
vertex sets of the $r$-partite graph $S$ 
contains exactly $h$ vertices, and that, \\
(2) for $I\in {\mathfrak{r}\choose [2]}$
 there is an injection from the index-$I$
visible colors of $S$ to the index-$I$ colors of ${\bf G}$.\\
(When the injection maps a visible color $\mathfrak{c}$ of $S$ to another
color $\mathfrak{c}'$ of ${\bf G}$, we simply write
$\mathfrak{c}=\mathfrak{c}'$ without presenting the injection
explicitly.) For $S\in {\mathcal S}_{h,{\bf G}}$, we denote by
${\mathbb V}_I(S)$ the set of index-$I$ visible edges. Write
${\mathbb V}_i(S):=\dot{\bigcup}_{I\in {\mathfrak{r}\choose i}}{\mathbb
V}_I(S)$ and ${\mathbb V}(S):={\mathbb V}_1(S)
\dot{\cup}{\mathbb V}_2(S).$
Clearly we have
\begin{eqnarray}
|\mathbb{V}_1(S)|\le rh \mbox{ and }|\mathbb{V}_2(S)|\le {r\choose 2}h^2.
\label{081021c}
\end{eqnarray}
\end{df}
\begin{df}[Partitionwise maps]
A {\bf partitionwise map} $\varphi: \dot{\bigcup}_i W_i\to \dot{\bigcup}_i
{\bf \Omega}_i$ is a map
 from $r$ disjoint
vertex sets $W_i,i\in \mathfrak{r},$ with $|W_i|<\infty$, to {\em the}
$r$ vertex sets (probability spaces) $ {\bf
\Omega}_i,i\in\mathfrak{r}$, such that each $w\in W_i$ is mapped
into ${\bf \Omega}_i.$ That is , any vertex is mapped to a vertex with the same index. 
We denote by
$\Phi((W_i)_{i\in\mathfrak{r}})$ or
$\Phi(\bigcup_{i\in\mathfrak{r}}W_i)$ the set of partitionwise
maps from $(W_i)_i.$ When $W_i=\{(i,1),\cdots, (i,h)\}$ or  when
$W_i$ are obvious and $|W_i|=h$, we denote it by $\Phi(h)$.
A partitionwise map is {\bf random} if and only if each $w\in W_i$ is 
independently mapped to a vertex in the probability space ${\bf
\Omega}_i$.
\end{df}
\begin{df}[Regularization]\label{070725}
Let $m\ge 0.$
Let ${\bf G}$ be a graph on ${\bf \Omega}$ and let $\varphi \in
\Phi(m).$ The {\bf regularization} of $\bf{G}$ by $\varphi$ is the
graph ${\bf G}/\varphi$ on ${\bf \Omega} $ obtained from ${\bf G}$
by redefining the color of each vertex ${\bf v}\in {\bf \Omega}_i$,
$i\in \mathfrak{r},$
by the $(1+(r-1)m)$-dimensional vector
\begin{eqnarray*}
\brkt{{\bf G}/\varphi}({\bf v}):=({\bf G}(
\{{\bf v}, {\bf u}\}) | 
{\bf u}={\bf v} \mbox{ or, }
{\bf u}\in 
{\bf \Omega}_j, j\in \mathfrak{r}\setminus\{i\},
\mbox{ is in the range of } \varphi ).
\end{eqnarray*}
\end{df}
Roughly speaking, the color of vertex ${\bf v}$ in ${\bf G}/\varphi$ 
is the information of the color-patterns of size-2 edges 
connecting the random vertex samplings and ${\bf v}$, together with the original color
${\bf G}({\bf v})$. (Here a size-2 edge means an edge which is not a single vertex but a pair of 
 vertices.)
\par
Note that edges of size $2$ (i.e. not vertices) do not get recolored in this process.
 Only vertices change their colors as the same as in the usual regularity lemma. 
\begin{df}[Regularity]\label{070724a}
Let ${\bf G}$ be a graph on $ {\bf \Omega}$. For
$\vec{\mathfrak{c}}=(\mathfrak{c}_J)_{J\subset I}\in {\rm
TC}_I({\bf G}), I\in {\mathfrak{r}\choose [2]}$, we define {\bf
relative density} by the conditional probability
\begin{eqnarray}
{\bf d}_{\bf G}(\vec{\mathfrak{c}}):= \Prob_{{\bf e}\in {\bf
\Omega}_I }[ {\bf G}({\bf e})=\mathfrak{c}_I | {\bf
G}(\partial{\bf e})= (\mathfrak{c}_J)_{J\subsetneq I} ].\label{081021d}
\end{eqnarray}
When $|I|=1$, in the above ${\bf e}$ is a vertex and the
conditional part is considered to always hold.
(Thus for $I=\{j\}$ and for an index-$\{j\}$ color $\mathfrak{c}_j\in {\rm TC}_1({\bf G})$, 
we have ${\bf d}_{\bf G}({\mathfrak{c}})
=\Prob_{{\bf v}\in {\bf \Omega}_j}[{\bf G}({\bf v})=\mathfrak{c}_j],
$ i.e. how much portion of the vertices in ${\bf \Omega}_j$ have color $\mathfrak{c}_j.$
For $I=\{i,j\}$ and for $\vec{\mathfrak{c}}=(\mathfrak{c}_{ij}; \mathfrak{c}_i,\mathfrak{c}_j)$, 
we have ${\bf d}_{\bf G}(\vec{\mathfrak{c}})=
\Prob_{{\bf v}_i{\bf v}_j\in {\bf
\Omega}_I }[ {\bf G}({\bf v}_i{\bf v}_j)=\mathfrak{c}_{ij}\ |\ {\bf
G}({\bf v}_i)= \mathfrak{c}_i \mbox{ and }
{\bf
G}({\bf v}_j)= \mathfrak{c}_j
].
$
)
\par
For a positive integer $h$ and $\epsilon\ge 0$, we say that
 ${\bf G}$ is
 {\bf $(\epsilon,h)$-regular}
if and only if there exists a function ${\delta}: {\rm TC}_2({\bf
G})\to [0,\infty)$
such that
\begin{eqnarray}
\hspace{-5mm} {\rm (i)}& \Prob_{\phi\in\Phi(h)}
[{\bf G}(\phi(e))=S(e)\,, \forall e\in {\mathbb V}(S)]\nonumber
\\
& =
\displaystyle\prod_{e\in {\mathbb V}_1(S)}  {\bf d}_{\bf G}(S\ang{e})
\displaystyle\prod_{e\in {\mathbb V}_2(S)} \brkt{ {\bf d}_{\bf G}(S\ang{e})
\dot{\pm} \delta(S\ang{e}) },& \forall S\in {\mathcal
S}_{h,{\bf G}}, \label{lj19}
\\
\hspace{-5mm} {\rm (ii)}& \Ex_{{\bf e}\in {\bf
\Omega}_I}[\delta({\bf G}\ang{\bf e})]\le \epsilon /|{\rm
C}_I({\bf G})|, & \forall I\in {\mathfrak{r}\choose 2},
\label{lj26}
\end{eqnarray}
where $a\dot{\pm}b$ denotes a suitable number $c$ satisfying
$\max\{0,a-b\}\le c\le \min\{1,a+b\}$.\\ Denote by ${\bf
reg}_{h}({\bf G})$ the minimum value of $\epsilon$ such that ${\bf
G}$ is $(\epsilon,h)$-regular.
\end{df}

{\bf Remark.} Roughly speaking, (i) measures how far from random
the graph {\bf G} is with respect to containing the expected
number of copies of the (colored) subgraphs $S \in
\mathcal{S}_{h,{\bf G}}$. The smaller $\delta$ is, the closer {\bf
G} is to being random. When $\delta\equiv 0$, then {\bf G} behaves
exactly like a random graph. On the other hand, if we take $\delta
\equiv 1$ then (i) is automatically satisfied. Condition (ii) places an
upper bound on the size of $\delta$.
Our proof will yield the main theorem even if we replace
the right-hand side of (ii) by $g_I(|{\rm C}_I({\bf G})|)$
for any fixed functions $g_I>0$, for example,
$g_I(x)=x^{-1/\epsilon}.$
\par
Our main theorem is as follows.
\begin{thm}[Main]\label{main}
For any $r\ge 2,\, h,\, \vec{b}=(b_1, b_2)$, and
$\epsilon>0,$ there exist an (increasing) function $
m:\naturalset\to\naturalset $
 and an integer
$\tilde{n}$ satisfying the following:
\par
If
 ${\bf G}$  is
a $\vec{b}$-colored ($r$-partite) graph on $ {\bf \Omega}$ then
\begin{eqnarray*}
\Ex_{n} \Ex_{\varphi} [{\bf reg}_{h}({\bf G}/{\varphi})] \le
\epsilon
\end{eqnarray*}
where $n$ is chosen randomly in $[0,\tilde{n}-1]$ and where
$\varphi \in \Phi(m(n))$ is random.
\end{thm}

Note that $m$ and $\tilde{n}$ depend only on $r,h,\vec{b}$ and
$\epsilon$ and are independent of everything else (including ${\bf
\Omega}$). Since $m$ is increasing, we put $\widetilde{m}:= m(\tilde{n})\ge m(n)$ and get:
\begin{cor}[Regularity Lemma]
\label{lj21} For any $r\ge 2,\,h,\,\vec{b}=(b_1,b_2)$, and
$\epsilon>0,$ there exists an integer $ \widetilde{m} $ such that
if
 ${\bf G}$  is
a $\vec{b}$-colored ($r$-partite) graph on $ {\bf \Omega}$ then
for some integer $m\le \widetilde{m}$, we have
\begin{eqnarray}
\Ex_{\varphi \in \Phi(m)} [{\bf reg}_{h}({\bf G}/\varphi)] \le
\epsilon.\label{081215}
\end{eqnarray}
In particular, when (\ref{081215}) holds, if we pick a map 
${\varphi}\in\Phi(m)$ randomly then with probability at least 
$1-\sqrt{\epsilon}$, we have ${\bf reg}_{h}({\bf G}/{\varphi})\le \sqrt{\epsilon},$ 
thus 
${\bf G}/{\varphi}$ is $(\sqrt{\epsilon},h)$-regular.
\end{cor}
It is important that the above integer $m$ is bounded by a constant $\widetilde{m}$ 
independent from ${\bf G}$ but the exact value of $m$ itself depends on ${\bf G}$.
Note that, in (cannonical) property testing, 
the exact value $m$ is also independent from ${\bf G}.$
This is a new critical idea which has never been previously while some 
had felt that property testing and graph regularization seem to have a close relation (ex. \cite{AFNS}).
\par
Of course, we can rewrite the above results for non-partite graphs.
\section{Proof of the Main Theorem}Before we proceed with the proof
of the Main Theorem, we will need to establish two lemmas. We
admit that they may appear a bit technical and unmotivated at this
point, but their use will be clearer once we see how they are used
in the main proof.
\subsection{Two lemmas and their proofs}
\begin{df}[Notation for the lemmas]\label{081229}
Let ${\bf G}$ be an ($r$-partite) graph on ${\bf \Omega}$. For two edges ${\bf
e}, {\bf e'}\in {\bf \Omega}_I,$ we abbreviate $ {\bf G}({\bf e})
= {\bf G}({\bf e'})$ and $ {\bf G}(\partial {\bf e}) = {\bf
G}(\partial {\bf e'}) $ by $ {\bf e}\stackrel{\bf G}{\approx}{\bf
e'}$ and $ {\bf e}\stackrel{\partial {\bf G}}{\approx}{\bf e'},$
respectively.
\par
An {\bf $h$-error function} of ${\bf G}$ is a function
$\delta:\bigcup_{I\in {\mathfrak{r}\choose 2}}{\rm TC}_I({\bf
G})\to [0,\infty)$ satisfying (\ref{lj19}) for all $S\in {\mathcal
S}_{h,{\bf G}}$.
\par
Denote by $\dbrkt{\cdots}$ the Iverson bracket, i.e., it equals
$1$ if the statement in the bracket holds, and $0$ otherwise.
\end{df}
\begin{lemma}[Correlation bounds counting error]\label{lj23c}
For any graph ${\bf G}$ on ${\bf \Omega}$ and for 
$S\in {\mathcal S}_{h,{\bf G}}$, we
have
\begin{eqnarray*}
&& \left| \Prob_{\phi\in\Phi(h)}\left[\left. {\bf
G}(\phi(e))=S(e), \,\forall e\in {\mathbb V}_2(S) \right| {\bf
G}(\phi(v))=S(v),\,\forall v\in {\mathbb V}_{1}(S)
\right]\nonumber - \prod_{e\in {\mathbb V}_2(S)} {\bf d}_{\bf
G}(S\ang{e}) \right|\\
&\le& |{\mathbb V}_2(S)| \max_{\emptyset\not=D\subset {\mathbb
V}_2(S)} \left| \Ex_{\phi\in\Phi(h)}\left[\left. \prod_{e\in D}
\brkt{ \dbrkt{{\bf G}(\phi(e))=S(e)} - {\bf d}_{\bf G}(S\ang{e}) }
\right| {\bf G}(\phi(v))=S(v),\,\forall v\in {\mathbb V}_{1}(S)
\right] \right| . \label{lj15b}
\end{eqnarray*}
\end{lemma}
\Proof [Tool: Nothing] We will prove this by induction on $|{\mathbb V}_2(S)|.$ If
$|{\mathbb V}_2(S)|\le 1$ then the statement is trivial, since in
this case, the expression on the left-hand side of the inequality
is $0$. So let us assume that $\left| {\mathbb V_2(S)}\right| \geq
2$ and that the result holds for all smaller values of $\left|
{\mathbb V_2(S)}\right|$. Let $d_e:={\bf d}_{\bf G}(S\ang{e})$,
and let $\eta$ be the maximum part of the desired right-hand side.
Then for $D:=\mathbb{V}_2(S)$ we have 
\begin{eqnarray*}
[-\eta,\eta]&\ni&
\Ex_{\phi\in\Phi(h )=\Phi(V(S))}\left[\left. \prod_{e\in {\mathbb
V}_2(S)}
\brkt{
\dbrkt{{\bf G}(\phi(e))=S(e)} - {\bf d}_{\bf G}(S\ang{e})
}
\right| {\bf G}(\phi(v))=S(v),\,\forall v\in {\mathbb V}_{1}(S)
\right]\nonumber\\
&=&  \Ex_{\phi\in\Phi(h )}\left[\left. \prod_{e\in {\mathbb
V}_2(S)} \dbrkt{{\bf G}(\phi(e))=S(e) } \right| {\bf
G}(\phi(v))=S(v),\,\forall v \in {\mathbb V}_{1}(S) \right]\nonumber
\\
&&+ \sum_{\emptyset\not=D\subset {\mathbb V}_2(S)}
\brkt{\prod_{e\in D} (-d_e) } \Ex_{\phi\in\Phi(h) }\left[ \left.
\prod_{e\in {\mathbb V}_2(S)\setminus D} \dbrkt{{\bf
G}(\phi(e))=S(e) } \right| {\bf G}(\phi(v))=S(v),\forall v \in
{\mathbb V}_{1}(S) \right], \label{081021}
\end{eqnarray*}
expanding the product and using the
linearity of expectation and the definition of $d_e$.
Now we will focus on second term above. Since the value of
$\dbrkt{{\bf G}(\phi(e))=S(e) }$ is $0$ or $1$, we can replace
$\Ex$ by $\Prob$, and consequently, apply the induction hypothesis
(since $D$ is nonempty). 
Consider a complex $S^-$ with $\mathbb{V}_2(S^-)=\mathbb{V}_2(S)\setminus D$ by 
invisualizing the edges in $D$ of $S$. 
\\[.1in]
Using the inductive hypothesis for complex $S^-$ 
 in the place of $S$, we rewrite the second term and obtain 
\begin{eqnarray*}
&& \Ex_{\phi\in\Phi(h )}\left[\left. \prod_{e\in {\mathbb V}_2(S)}
\dbrkt{{\bf G}(\phi(e))=S(e) } \right| {\bf
G}(\phi(v))=S(v),\,\forall v\in {\mathbb V}_{1}(S) \right]
\nonumber
\\
&
\stackrel{\rm I.H.}
{=}
& - \sum_{\emptyset\not=D\subset {\mathbb
V}_2(S)} \brkt{\prod_{e\in D} (-d_e) } \brkt{ 
\brkt{\prod_{e\in {\mathbb
V}_2(S)\setminus D} d_e}\pm
|{\mathbb V}_2(S^-)|
 \eta
 } \quad 
\dot{\pm}\quad \eta
\nonumber
\\
&=& 
- \brkt{\brkt{ \prod_{e\in {\mathbb V}_2(S)} d_e}\pm
|{\mathbb V}_2(S^-)| \eta
} \sum_{\emptyset\not=D\subset {\mathbb V}_2(S)} \brkt{\prod_{e\in
D} (-1) } \quad \dot{\pm}\, \eta \quad (
\because |d_e|\le 1)
\nonumber
\\
&=& 
- \brkt{ 
\brkt{\prod_{e\in {\mathbb V}_2(S)} d_e}\pm
(|{\mathbb V}_2(S)|-1) \eta
}
\brkt{ (1-1)^{|{\mathbb V}_2(S)|} -1 } 
\quad \dot{\pm}\,\eta 
 \quad (\because 
|\mathbb{V}_2(S)|>|\mathbb{V}_2(S^-)|
)
\nonumber
\\
&=&
\brkt{
\prod_{e\in {\mathbb V}_2(S)}
d_e}\pm|{\mathbb V}_2(S)|\eta.
\end{eqnarray*}
\lemmaqed
We will use the following form of the Cauchy-Schwarz.
\begin{fact}[Cauchy-Schwarz inequality]\label{080827a}
For a random variable $X$ on a probability space ${\Omega}$ if 
an equivalent relation $\approx$ on $\Omega$ is a refinement of another equivalent relation $\sim$ 
on $\Omega$ then
\begin{eqnarray}
\Ex_{{\omega_0}\in {\Omega}}\brkt{\Ex_{\omega\in \Omega}
[X({\omega})|{\omega}\approx {\omega_0}]}^2
\ge 
\Ex_{\omega_0\in\Omega}\brkt{\Ex_{\omega\in\Omega}[X({\omega})| {\omega}\sim {\omega_0}]
}^2.\label{080827}
\end{eqnarray}
\end{fact}
\Proof 
By the Cauchy-Schwarz 
(i.e. $\Ex[X^2]\Ex[Y^2]\ge (\Ex[XY])^2$), 
we have 
$
\Ex_{\omega_0}\brkt{\Ex_{\omega}[X({\omega})|{\omega}\approx {\omega_0}]}^2
=
\Ex_{\omega_0}\left[
\Ex_{\omega'}\left[\left.
\brkt{\Ex_{\omega}[X({\omega})|{\omega}\approx {\omega'}]
}^2\right|{\omega'}\sim {\omega_0}
\right]\right]
=
\Ex_{\omega_0}\left[
\Ex_{\omega'}[1^2|{\omega'}\sim {\omega_0}]
\cdot
\Ex_{\omega'}\left[\left.
\brkt{\Ex_{\omega}[X({\omega})|{\omega}\approx {\omega'}]
}^2\right|{\omega'}\sim {\omega_0}
\right]\right]
\stackrel{CS}{\ge}
\Ex_{\omega_0}\left(
\Ex_{\omega'}[
{1\cdot\Ex_{\omega}[X({\omega})|{\omega}\approx {\omega'}]
}|{\omega'}\sim {\omega_0}
]\right)^2= 
\Ex_{\omega_0}\brkt{\Ex_{\omega}[X({\omega})| {\omega}\sim {\omega_0}]
}^2.
$
\factqed
With this fact and Definition \ref{081229}, we next tackle
\begin{lemma}[Mean square bounds correlation]
\label{alj08} 
Let $h$ and $m$ be positive integers and ${\bf G}$
an $r$-partite graph on the vertex set ${\bf \Omega}$. Let $S\in
{\mathcal S}_{h,{\bf G}}$ and let $F_e:{\rm C}_I({\bf G}) \to
[-1,1]$ be a function for each $I\in {\mathfrak{r}\choose 2}$ and
for each $e\in {\mathbb V}_I(S).$ For any $I\in {\mathfrak{r}
\choose 2}$ and $e_0\in {\mathbb V}_I(S),$ we have
\begin{eqnarray}
&& \brkt{ \Ex_{\phi\in\Phi(h)}\left[ \prod_{e\in {\mathbb V}_2(S)}
F_e\brkt{\phi(e)} \prod_{v\in {\mathbb V}_{1}(S)} \llbracket {\bf
G}(\phi(v))=S(v) \rrbracket \right]}^2\nonumber
\\
&\le& \Ex_{\varphi\in \Phi(mh)} \Ex_{{\bf e^*}\in {\bf \Omega}_I }
[ \brkt{ \Ex_{{\bf e}\in {\bf \Omega}_I }[F_{e_0}\brkt{\bf e}
\llbracket {\bf G}(\partial {\bf e})=S(\partial e_0)\rrbracket |
{\bf e}\stackrel{
\partial {\bf G}/\varphi}{\approx}{\bf e^*}]
}^2 ]\nonumber
\\
&&\cdot \brkt{ \prod_{v\in {\mathbb V}_{1}(S)} {\bf d}_{\bf
G}^{}(S\ang{v}) } \brkt{ \brkt{ \prod_{v\in {\mathbb
V}_{1}(S),v\not\in e_0} {\bf d}_{\bf G}^{}(S\ang{v}) } + { 1 \over
m} }
\label{070907d}
\end{eqnarray}
where $\phi, \varphi$ are random and where we abbreviate $F_e({\bf
G}({\bf e}))$ by $F_e({\bf e})$.
\par
In particular, if we suppose $ {1\over m} \le \prod_{v\in {\mathbb
V}_{1}(S), v\not\in e_0 } {\bf d}_{\bf G}(S\ang{v}) $ (i.e., $m$
is large) then
\begin{eqnarray}
&&
\brkt{
\Ex_{\phi\in\Phi(h)}\left[
\left.
\prod_{e\in {\mathbb V}_2(S)}
F_e\brkt{\phi(e)}
\right|
{\bf G}(\phi(v))=S(v) \, \forall
v\in {\mathbb V}_{1}(S)
\right]
}^2\nonumber
\\
&\le&
2
\Ex_{\varphi\in \Phi(mh)}
\Ex_{{\bf e^*}\in {\bf \Omega}_I
}
[
\brkt{
\Ex_{{\bf e}\in {\bf \Omega}_I
}[F_{e_0}\brkt{\bf e}
|
{\bf e}\stackrel{
\partial {\bf G}/\varphi}{\approx}{\bf e^*}]
}^2
|{\bf G}(\partial {\bf e^*})=S(\partial e_0)
].
\label{070907c}
\end{eqnarray}
\end{lemma}
\Proof [Tools: Cauchy-Schwarz, Fact \ref{080827a}]
Fix $I_0\in {\mathfrak{r}\choose 2}$ and $e_0\in {\mathbb
V}_{I_0}(S)$. For $\phi\in\Phi(V(S)\setminus e_0)$ and for ${\bf
e_0}\in {\bf \Omega}_{I_0}$, we define the (extended) function
$\phi^{({\bf
e_0})}\in\Phi(V(S))$ such that:\\
 (i) each $v\in e_0$ is mapped to
the corresponding ${\bf v}\in {\bf e_0}$ with the index of $v$,
and
 that,\\
  (ii) each $v\in V(S)\setminus e_0$ is mapped to
$\phi(v)$.\\ 
(That is, when we have a map $\phi$ from all but two vertices $e_0$, 
we extend it by assigning two vertices $e_0$ to two vertices ${\bf e_0}.$
) 
For an $m$-tuple of maps $\vec{\varphi}=(\varphi_i)_{i\in [m]}$ with
$\varphi_i\in\Phi(V(S)\setminus e_0),$ we define an equivalence
relation $\stackrel{\vec{\varphi}}{\sim}$ on ${\bf \Omega}_{I_0}$
by the condition that
\begin{eqnarray}
{\bf e}\stackrel{\vec{\varphi}}{\sim}{\bf e'}
\mbox{ if and only if } \varphi_i^{({\bf e})}(e) \stackrel{\bf G}{\approx}
\varphi_i^{({\bf e'})}(e)  \quad  \forall e\in
\mathbb{V}(S)\setminus\{e_0\},
\forall i\in [m].\label{lj23}
\end{eqnarray}
(Note that $V(S)\setminus e_0$ is a vertex set while $\mathbb{V}(S)\setminus \{e_0\}$ is 
an edge set.
Since the right-hand side clearly holds for $e$ with $e\cap e_0=\emptyset,$
it is enough to check only the edges $e\in \mathbb{V}(S)$
with $|e\cap e_0|=1.$)\\
Let $S^{(1)},\cdots,S^{(m)}$ and $e_0^{(1)},\cdots,e_0^{(m)}$ be
$m$ copies of $S$ and $e_0$. 
For $\vec{\varphi}=(\varphi_i)_{i\in
[m]}$ with $\varphi_i\in\Phi (V(S^{(i)})\setminus e_0^{(i)})$, let
$\varphi^*\in
\Phi(mh)=\Phi(V(S^{(1)})\dot{\cup}\cdots\dot{\cup}V(S^{(m)}))$ be
an extended function of $\varphi_i$'s, i.e.,
$\varphi^*(v)=\varphi_i(v)$ for all $v\in V(S^{(i)})\setminus
e_0^{(i)}, i\in [m]$. 
Then it is not hard to see that
\begin{eqnarray}
{\bf e}\stackrel{\partial {\bf G}/\varphi^*}{\approx}{\bf e'}
\mbox{ implies } {\bf e}\stackrel{\vec{\varphi}}{\sim}{\bf e'}.
\label{lj23a}
\end{eqnarray}
(To see this, observe that 
if $I_0=\{1,2\}$ and 
$\{{\bf v}_1, {\bf v}_2\}
\stackrel{\partial {\bf G}/\varphi^*}{\approx}
\{{\bf v'}_1, {\bf v'}_2\},
$ or equivalently 
${\bf v}_j\stackrel{{\bf G}/\varphi^*}{\approx}{\bf v'}_j (j=1,2),
$ then $\{{\bf v}_j, \varphi^*(v)\}\stackrel{{\bf G}}{\approx}\{{\bf v'}_j, \varphi^*(v)\}
$ for all $v\in 
V(S^{(i)})$ having no index $j$ (i.e. $v\in V_{j'}(S^{(i)}), j'
\not=j$).
Since $\varphi^*(v)=\varphi_i(v)$ if $v\not\in e_0^{(i)}$, 
we see $\varphi_i^{(\{{\bf v}_1,{\bf v}_2\})}(e)\stackrel{\bf G}{\approx}
\varphi_i^{(\{{\bf v'}_1,{\bf v'}_2\})}(e)
$ for all $e\in \mathbb{V}_2(S)\setminus \{e_0\}$ with $|e\cap e_0|=1,$ 
implying (\ref{lj23a}) by (\ref{lj23}).)
\par
Let $F^*_{e_0}({\bf e}):= F_{e_0}({\bf e})
\dbrkt{{\bf G}(\partial{\bf e})=S(\partial e_0)}
$ and let
\begin{eqnarray*}
F^*(\phi):= \prod_{e\in {\mathbb V}_2(S)\setminus \{e_0\}}
F_e\brkt{\phi(e)}\prod_{v\in {\mathbb V}_1 (S)} \dbrkt{{\bf
G}(\phi(v))=S(v)}.
\end{eqnarray*}
Then, since $\dbrkt{\cdots}^2=\dbrkt{\cdots},$ 
the left-hand side of (\ref{070907d}) becomes 
\begin{eqnarray*}
&& \brkt{ \Ex_{\phi\in\Phi(h)}\left[ F_{e_0}(\phi(e_0))
\prod_{e\in {\mathbb V}_2(S)\setminus\{e_0\}}
F_e\brkt{\phi(e)} \cdot 
\dbrkt{{\bf G}(\partial(\phi(e_0)))=S(\partial e_0)}^2
\prod_{v\in {\mathbb V}_{1}(S):v\not\in e_0} \llbracket {\bf
G}(\phi(v))=S(v) \rrbracket \right]}^2
\\ &=& 
\brkt{ \Ex_{\phi\in\Phi(h)}\left[
F^*_{e_0}\brkt{\phi(e_0)}F^*(\phi) \right]}^2\nonumber
\mbox{\hspace{1.5in}(by the definition of $F_{e_0}^*$ and $F^*$})\\
&=& \brkt{ \Ex_{{\bf e_0}\in {\bf \Omega}_{I_0},
\phi\in\Phi(V(S)\setminus e_0)}\left[ F^*_{e_0}\brkt{\bf
e_0}F^*(\phi^{({\bf e_0})}) \right] }^2
\nonumber\\
&& \mbox{\hspace{.2in} 
(since the two expectations are taken 
over random choices of $2+(rh-2)=rh$ vertices in $V(S)$)}
\\
&=& \brkt{\Ex_{\vec{\varphi}=(\varphi_i)_{i\in [m]}\in
(\Phi(V(S)\setminus e_0))^m} \Ex_{{\bf e_0}\in {\bf \Omega}_{I_0}}
\left[ F^*_{e_0}\brkt{\bf e_0} \Ex_{i\in [m]}[F^*(\varphi_i^{({\bf
e_0})})] \right] }^2
\nonumber\\
\nonumber\\
&=& \brkt{\Ex_{\vec{\varphi}=(\varphi_i)_{i\in [m]}\in
(\Phi(V(S)\setminus e_0))^m} \Ex_{{\bf e_0}\in {\bf \Omega}_{I_0}}
\left[ \Ex_{{\bf e}\in {\bf \Omega}_{I_0}} \left[
F^*_{e_0}\brkt{\bf e} \Ex_{i\in [m]}[F^*(\varphi_i^{({\bf e})})]
\left|{\bf e}\stackrel{\vec{\varphi}}{\sim}{\bf e_0}\right.
\right] \right] }^2\\
&& \mbox{\hspace{.2in} (since for any random variable $X$ and the equivalence classes 
$C_i$ by $\stackrel{\vec{\varphi}}{\sim},$
}\\
&& \mbox{\hspace{.2in} 
$\Ex_{\bf e_0}\Ex_{\bf e}[X({\bf e})|{\bf e}\stackrel{\vec{\varphi}}{\sim}{\bf e_0}]
=\sum_i\Prob_{\bf e_0}[{\bf e_0}\in C_i]
\Ex_{\bf e}[X({\bf e})|{\bf e}\in C_i]
=\Ex_{\bf e_0}[X({\bf e_0})]
$
)}\\
\nonumber &
\stackrel{(\ref{lj23})}{=}
& \brkt{\Ex_{\vec{\varphi}=(\varphi_i)_{i\in [m]}\in
(\Phi(V(S)\setminus e_0))^m} \Ex_{{\bf e_0}\in {\bf \Omega}_{I_0}}
\left[ \Ex_{{\bf e}\in {\bf \Omega}_{I_0}} [ F^*_{e_0}\brkt{\bf e}
|{\bf e}\stackrel{\vec{\varphi}}{\sim}{\bf e_0} ] \Ex_{i\in
[m]}[F^*(\varphi_i^{({\bf e_0})})] \right] }^2\\ &&
\mbox{\hspace{1.5in} (since $F^*(\varphi^{({\bf
e})}_i)=F^*(\varphi^{({\bf e_0})}_i)$ when ${\bf
e}\stackrel{\vec{\varphi}}{\sim}{\bf e_0}$)}\\
\nonumber &\leq& { \Ex_{\vec{\varphi} } \Ex_{{\bf e_0}} \left[
\brkt{ \Ex_{{\bf e}\in {\bf \Omega}_{I_0}} [ F^*_{e_0}\brkt{\bf e}
| {\bf e}\stackrel{\vec{\varphi}}{\sim}{\bf e_0} ] }^2 \right]
\cdot \Ex_{\vec{\varphi}=(\varphi_i)_i } \Ex_{{\bf e_0}} \left[
\brkt{ \Ex_{i\in [m]}[ F^*(\varphi_i^{({\bf e_0})}) ] }^2 \right]
}
\mbox{ (by \underline{Cauchy-Schwarz})}\\
&
\stackrel{(\ref{lj23a}),(\ref{080827})}{\leq}& \Ex_{\varphi^*\in\Phi(mh) } \Ex_{{\bf e_0}} \left[ \brkt{
\Ex_{{\bf e}\in {\bf \Omega}_{I_0}} [ F^*_{e_0}({\bf e}) | {\bf
e}\stackrel{\partial {\bf G}/\varphi^*}{\approx}{\bf e_0} ] }^2
\right] \cdot \Ex_{\vec{\varphi}=(\varphi_i)_i }\Ex_{{\bf e_0}\in
{\bf \Omega}_{I_0}} \left[ \Ex_{i,j\in [m]}[ F^*(\varphi_i^{({\bf
e_0})}) F^*(\varphi_j^{({\bf e_0})}) ] \right].
\end{eqnarray*}
So, the first term in this last line is the first term in our
desired inequality. We now focus on the second term. 
Since $|F^*(\cdot)|\le 1,$ we have  
\begin{eqnarray*}
&& \Ex_{\vec{\varphi}=(\varphi_i)_i }\Ex_{{\bf e_0}\in {\bf
\Omega}_{I_0}} \left[ \Ex_{i,j\in [m]}[ F^*(\varphi_i^{({\bf
e_0})}) F^*(\varphi_j^{({\bf e_0})}) ] \right]\\
&=& \Ex_{{\bf e_0}\in {\bf \Omega}_{I_0}} \left[{m-1\over m}
\Ex_{i\not=j\in [m]} \Ex_{\vec{\varphi}=(\varphi_i)_i} [
F^*(\varphi_i^{({\bf e_0})}) F^*(\varphi_j^{({\bf e_0})}) ]
+{1\over m} \Ex_{i\in [m]} \Ex_{\vec{\varphi}=(\varphi_i)_i}[
\brkt{ F^*(\varphi_i^{({\bf e_0})}) }^2 ] \right]\\
&\leq& \Ex_{{\bf e_0}\in {\bf \Omega}_{I_0}}
\Ex_{\varphi_1,\varphi_2\in \Phi(V(S)\setminus e_0)} [\left|
F^*(\varphi_1^{({\bf e_0})}) F^*(\varphi_2^{({\bf e_0})})\right| ]
+{1\over m}\Ex_{{\bf e_0}\in {\bf \Omega}_{I_0}} \Ex_{\varphi_1\in
\Phi(V(S)\setminus e_0)} [ \left|F^*(\varphi_1^{({\bf
e_0})})\right| ]\\
&\leq& \Ex_{{\bf e_0}\in {\bf \Omega}_{I_0}}
\Ex_{\varphi_1,\varphi_2\in \Phi(V(S)\setminus e_0)} \left[
\prod_{v\in {\mathbb V}_{1}(S)} \left| \dbrkt{{\bf
G}(\varphi_1^{({\bf e_0})}(v))=S(v)} \dbrkt{{\bf
G}(\varphi_2^{({\bf e_0})}(v))=S(v)} \right| \right]\\
&& +{1\over m} \Ex_{\varphi\in \Phi(V(S))} \left[ \prod_{v\in
{\mathbb V}_{1}(S)} \dbrkt{{\bf G}(\varphi(v))=S(v)} \right].
\quad\quad\quad (\mbox{by the definition of $F^*$ since }|F_e|\le
1)
\end{eqnarray*}
Looking at the second term first, this can be written as
\begin{eqnarray*}
&& {1\over m}
\prod_{v\in {\mathbb V}_1(S)}
\Prob_{\varphi \in \Phi({\mathbb V(S)})}  \left[{\bf
 G}(\varphi(v))=S(v)\right]\quad \mbox{(since $\varphi$ maps all $v \in
\Phi({\mathbb V_{1}(S)})$
 independently)} \\
 &=& {1\over m}\prod_{v \in {\mathbb V}_{1}(S)} {\bf d}_{{\bf
 G}}(S\ang{v})\quad \mbox{(by the definition of ${\bf d}_{{\bf
 G}}(S\ang{v})$
)
}.\\
\end{eqnarray*}
In a similar way, we can interpret the first term as computing the
probability that $2+2(rh-2)$ random (visible or invisible) vertices chosen independently 
will have vertex colors in ${\bf G}$ which match
those of their corresponding vertices in $S$. This
probability can be written as
\begin{eqnarray*}
\prod_{v\in {\mathbb V}_{1}(S), v\not\in e_0 } \Prob_{{\bf v} \in
\Omega_I} \left[ {\bf G}({\bf v}) = S(v) \right]^2 \prod_{v
\in e_0} \Prob_{{\bf v} \in
\Omega_I} \left[ {\bf G}({\bf v}) = S(v) \right]\\
= \prod_{v\in {\mathbb V}_{1}(S), v\not\in e_0} {\bf d}_{{\bf
G}}(S\ang{v})^2\,\prod_{v \in e_0} {\bf d}_{{\bf
 G}}(S\ang{v}).\\
\end{eqnarray*}
Putting all these observations together, the proof of the first part of
Lemma
\ref{alj08} is complete.
\par
Next, we show the last sentence of the lemma.
The left hand side of (\ref{070907c}) is at most
\begin{eqnarray*}
&&
\brkt{
\Ex_{\phi\in\Phi(h)}\left[
\prod_{e\in {\mathbb V}_2(S)}
F_e\brkt{\phi(e)}
\prod_{v\in {\mathbb V}_1(S)}
\llbracket
{\bf G}(\phi(v))=S(v)
\rrbracket
\right]
/
\Prob_{\phi\in\Phi(h)}\left[
{\bf G}(\phi(v))=S(v) \, \forall
v\in {\mathbb V}_{1}(S)
\right]
}^2\nonumber
\\
&=&\brkt{
\Ex_{\phi\in\Phi(h)}\left[
\prod_{e\in {\mathbb V}_2(S)}
F_e\brkt{\phi(e)}
\prod_{v\in {\mathbb V}_1(S)}
\llbracket
{\bf G}(\phi(v))=S(v)
\rrbracket
\right]
/
\brkt{
\prod_{v\in {\mathbb V}_{1}(S)}
{\bf d}_{\bf G}(S\ang{v})
}}^2
\\
&\stackrel{(\ref{070907d})}{\le }&
\Ex_{\varphi\in \Phi(mh)}
\Ex_{{\bf e^*}\in {\bf \Omega}_I
}
[
\brkt{
\Ex_{{\bf e}\in {\bf \Omega}_I
}[F_{e_0}\brkt{\bf e}
\llbracket
{\bf G}(\partial {\bf e})=S(\partial e_0)\rrbracket
|
{\bf e}\stackrel{
\partial {\bf G}/\varphi}{\approx}{\bf e^*}]
}^2
|{\bf G}(\partial {\bf e^*})=S(\partial e_0)
]
\\
&&
\cdot
\brkt{
\prod_{v\in {\mathbb V}_{1}(S), v\in e_0}
{\bf d}_{\bf G}(S\ang{e})
}
\brkt{
\prod_{v\in {\mathbb V}_{1}(S)}
{\bf d}_{\bf G}(S\ang{v})
}
\brkt{
\brkt{
\prod_{v\in {\mathbb V}_{1}(S),v\not\in e_0}
{\bf d}_{\bf G}(S\ang{v})
}
+
{
1 \over m}
}
\\
&&
/
\brkt{
\prod_{v\in {\mathbb V}_{1}(S)}
{\bf d}_{\bf G}(S\ang{v})
}^2
\mbox{ 
($\because$ ${\bf e}\stackrel{\partial {\bf G}/\varphi}{\approx}{\bf
e^*}$ and ${\bf G}(\partial {\bf e^*})=S(\partial e_0)$ imply
${\bf G}(\partial {\bf e})=S(\partial e_0)$)
}
\\
&\le &
\Ex_{\varphi\in \Phi(mh)}
\Ex_{{\bf e^*}\in {\bf \Omega}_I
}
[
\brkt{
\Ex_{{\bf e}\in {\bf \Omega}_I
}[F_{e_0}\brkt{\bf e}
|
{\bf e}\stackrel{
\partial {\bf G}/\varphi}{\approx}{\bf e^*}]
}^2
|{\bf G}(\partial {\bf e^*})=S(\partial e_0)
]
\\
&&
\cdot
\brkt{
1
+
{
1 \over m}
\brkt{
\prod_{v\in {\mathbb V}_{1}(S),v\not\in e_0}
{\bf d}_{\bf G}(S\ang{v})
}^{-1}
}.
\nonumber
\end{eqnarray*}
The assumption on $m$
now completes the proof of (\ref{070907c}). \lemmaqed
\subsection{The body of our proof}

\begin{df}[Notation for the proof]
Write $c_i({\bf G}):=\max_{I\in {{\mathfrak r}\choose i}} |{\rm
C}_I({\bf G})|$ for $i=1,2$. For $\vec{b}=(b_1,b_2)$ and
an integer $m$, we write
$\vec{B}(\vec{b},m):=({B}_i(\vec{b},m))_{i\in [2]}$ where
$ B_1(\vec{b},m):=b_1 \cdot b_{2}^{(r-1)m}$ and $B_2(\vec{b},m)
:= b_2$. Recalling the definition of regularization ${\bf G}/\varphi$, 
it is easy to see that if ${\bf G}$ is a $\vec{b}$-colored graph then
\begin{eqnarray}
c_i({\bf G}/\varphi) \le B_i(\vec{b},m),\quad \forall i=1,2,\,
\forall \varphi\in\Phi(m).\label{lj24k}
\end{eqnarray}
(If $i=2$, it is obvious since regularization does not recolor any size-2 edge. 
If $i=1,$ the new color of a vertex is determined by its original color and by the colors 
of the edges 
connecting the vertex and the $(r-1)m$ random vertices.
)
\end{df}
Suppose we are given some fixed $h \geq 1$, $\epsilon > 0$ and
$\vec{b}$. Our job will be to define suitable functions $m$ and
$\delta$, and a suitable integer $\tilde{n}$, so that (\ref{lj19})
and (\ref{lj26}) are satisfied. This we now do.
\par
\par
$\bullet$ {\bf [Definition of the sample-size functions] } Set
$m_{h,\vec{b},\epsilon}(0) := m(0):=0 $. Define
$\tilde{n}_{2,h,\vec{b},\epsilon}=\tilde{n}$ to be large enough so
that
\begin{eqnarray}
C\,b_2 \sqrt{ b_2\over \tilde{n} }\le {\epsilon \over 2{r\choose
2}} \label{lj24g}
\end{eqnarray}
where
\begin{eqnarray}
C:=\sqrt{2} {r\choose 2}h^2 \brkt{b_2\over \sqrt{\epsilon_1}
}^{{r\choose 2}h^2-1} \mbox{ and } \epsilon_1:=\brkt{\epsilon\over
 6\cdot b_2{r\choose 2}}^{2}.\label{lj24b}
\end{eqnarray}
(These expressions will appear in (\ref{071031a}) and (\ref{lj24c}).)
\par
We will define the function $m$ recursively as follows. Suppose
that $m(n)$ has been defined for some value of $n \geq 0$. Let
\begin{eqnarray}
M&:=&  \brkt{b_1b_2^{(r-1) m( n )} \over \sqrt{\epsilon_1}}^{rh}.
\label{lj24h}
\end{eqnarray}
(We will use the form (\ref{lj24h}) only once in (\ref{lj22b}).)
Define $m(n+1)$ so that
\begin{eqnarray}
m(n+1)&\ge& m( n ) +Mh = m( n ) +  \brkt{b_1b_2^{(r-1) m( n )} \over
\sqrt{\epsilon_1}}^{rh}h . \label{lj24a}
\end{eqnarray}\\[.1in]
Next, we define the error function $\delta$.
\\[.1in] $\bullet$ {\bf
[Definition of the error function]} For $\varphi\in\Phi(m(n))$, we
write ${\bf G^*}:={\bf G}/\varphi$ and we define the error
function $ \delta=\delta_{h,\epsilon,{\bf G^*}}$ inductively as
follows.
\par
First, define 
\begin{eqnarray}
\delta(\vec{\mathfrak{c}}):=0 \mbox{ and }
\eta(\vec{\mathfrak{c}}):=0 \mbox{ for all } \vec{\mathfrak{c}}\in {\rm
TC}_I({\bf G^*}) \mbox{ with } I\in {\mathfrak{r}\choose 1}=\mathfrak{r}.
\label{081021b}
\end{eqnarray}
Before defining $\delta(\vec{\mathfrak{c}})$ and
$\eta(\vec{\mathfrak{c}})$ for $\vec{\mathfrak{c}} \in {\rm
TC}_2({\bf G^*})$, we define \lq bad colors\rq\ ${\rm BAD}\subset
{\rm TC}({\bf G^*}).$ For $I\in {\mathfrak{r}\choose [2]}$, we
define ${\rm BAD}_I$ by the relation that
$\vec{\mathfrak{c}}=(\mathfrak{c}_J)_{J\subset I}\in {\rm BAD}_I$
if and only if
\begin{eqnarray}
\begin{array}{cccl}
{\bf d}_{\bf G^*}((\mathfrak{c}_J)_{J\subset I^*} )&\le &
\sqrt{\epsilon_1}/|{\rm C}_{I^*}({\bf G^*})|& \mbox{ for some } I^* 
\mbox{ with }
\emptyset\not=I^*\subset I.
\end{array}\label{lj24m}
\end{eqnarray}
Define ${\rm BAD}:=\bigcup_{I\in {\mathfrak{r}\choose [2]}}{\rm
BAD}_I$. A {\em bad edge} will mean a visible edge whose color is bad.
\par
For $ \vec{\mathfrak{c}}=(\mathfrak{c}_J)_{J\subset I} \in {\rm
TC}_2({\bf G^*})$, we define, using $M$ and $C$ of
(\ref{lj24b}) and (\ref{lj24h}), 
\begin{eqnarray}
\eta(\vec{\mathfrak{c}})&:=& \Ex_{\varphi'\in\Phi(Mh)}\Ex_{{\bf
e^*}\in {\bf \Omega}_I} [ \brkt{\Prob_{{\bf e}\in {\bf \Omega}_I}[
{\bf G^*}({\bf e})=\mathfrak{c}_I | {\bf e}\stackrel{\partial {\bf
G^*}/\varphi'}{\approx}{\bf e^*} ] - {\bf d}_{\bf
G^*}(\vec{\mathfrak c})}^2 |\, {\bf G^*}(\partial {\bf e^*})=
(\mathfrak{c}_J)_{J\subsetneq I} ],\label{061220}
\\
\delta(\vec{\mathfrak{c}})&:=& \left\{
\begin{array}{cl}
1 & \mbox{ if } \vec{\mathfrak{c}}\in{\rm BAD}_I ,
\\
C\sqrt{\eta(\vec{\mathfrak{c}})}, & \mbox{otherwise.}\label{lj25}
\end{array}
\right.
\end{eqnarray}\\[.1in]
First, we show that with the above specified choices for $m, \tilde{n}$
and $\delta$, (\ref{lj19}) is satisfied.\\
$\bullet$ {\bf [The qualification as an error function]} 
Clearly it is enough to show that
\begin{eqnarray}
&&
\Prob_{\phi\in\Phi(h)}
[{\bf G^*}(\phi(e))=S(e)\,, \forall e\in {\mathbb V}(S)]
\nonumber
\\
 &=&
\displaystyle\prod_{e\in {\mathbb V}_1(S)}  {\bf d}_{\bf G^*}(S\ang{e})
\displaystyle\prod_{e\in {\mathbb V}_2(S)} \brkt{ {\bf d}_{\bf G^*}(S\ang{e})
\dot{\pm} \delta(S\ang{e}) }
\label{lj24}
\end{eqnarray}
for any $S\in {\mathcal S}_{h,{\bf G^*}}.$
Furthermore without loss of generality, we can assume that
\begin{eqnarray}
S\ang{e}\not\in{\rm BAD} \mbox{ for any } e\in {\mathbb V}(S).
\label{071105}
\end{eqnarray}
(Indeed, we can show this by the induction on the number of bad edges in $S$.
Let a complex $S$ be given where $S$ contains a bad edge $e^*$. 
Firstly we suppose that there exist no bad vertices and thus 
 $e^*$ contains two different vertices (which are not bad).
By the induction hypothesis, 
(\ref{lj24}) holds for the complex $S^*$ 
obtained from $S$ by recoloring $e^*$ in the invisible color. 
Equality (\ref{lj24}) means that 
the real number the left hand side suggests belongs to the 
interval which the right-hand side suggests. Denote by $[p^-, p^+]$ this interval.
Again we reconstruct $S$ from $S^*$ by recoloring some invisible edges in the original 
bad color.
By this process from $S$ to $S^*$, 
the left hand side of (\ref{lj24}) will not increase 
(probably decrease because of added visible edges $e^*$)
and 
 the right-hand side will suggest interval $[0,p^+]$ because 
 ${\bf d}_{\bf G^*}(S\ang{e^*}
)\dot{\pm}\delta(S\ang{e^*})=[0,1]$ 
 by 
(\ref{lj25}).  
Then (\ref{lj24}) holds also for $S$. 
Secondly we suppose that the $e^*$ consists of a single bad vertex $v^*$. 
Then we recolor not only $v^*$ but also all edges containing $v^*$ in the 
invisible color. The same argument can be applied.
)
\\[.1in]
Fix such an $S\in {\mathcal S}_{h,{\bf G^*}}$. 
For any $e\in
{\mathbb V}_J(S)$ with $J\subset \mathfrak{r}$, it follows from 
(\ref{071105}), (\ref{lj24m}) and (\ref{081021b}) that
\begin{eqnarray}
{\bf d}^{}_{\bf G^*}(S\ang{e})
> {\sqrt{\epsilon_1}\over
|{\rm C}_J({\bf G^*})|}>0 \quad(\mbox{if }|J|\le 2) \mbox{ and }
\delta(S\ang{e})=0 \,(\mbox{if }|J|=1).\label{lj22}
\end{eqnarray}\\[.1in]
Using (\ref{lj24k}), (\ref{lj24h}) and (\ref{lj22}), a 
straightforward computation gives
\begin{eqnarray}
{1\over M} \stackrel{(\ref{lj24k}),(\ref{lj24h})}{\le}
 \brkt{\sqrt{\epsilon_1}\over c_1({\bf G^*})
}^{|{\mathbb V}_1(S)|} 
\stackrel{(\ref{lj22})}{\le}
 \prod_{v\in {\mathbb V}_{1}(S)} {\bf
d}_{\bf G^*}(S\ang{v}) \leq \prod_{v\in {\mathbb V}_{1}(S),
v\not\in e_0 } {\bf d}_{\bf G^*}(S\ang{v})\label{lj22b}
\end{eqnarray}
for any $e_0\in\mathbb{V}_2(S).$
For any choice of $\emptyset \neq D \subset \mathbb{V}_2(S)$, we
 define $S'\in {\mathcal S}_{h,{\bf G^*}}$ so that
$\mathbb{V}_2(S')=D$ and $S'(e) = S(e) \, \forall e \in D$ and that $S'(v)=S(v) \ \forall 
v\in 
\mathbb{V}_1(S)=\mathbb{V}_1(S').$ 
Now,
applying Lemma \ref{alj08} for $S'$ with $ F_e({\bf e}):=
\dbrkt{{\bf G^*}({\bf e})=S(e) }-{\bf d}_{\bf G^*}(S\ang{e})$, we have 
\begin{eqnarray}
&& \brkt{ \Ex_{\phi\in\Phi(h)}\left[\left. \prod_{e\in
\mathbb{V}_2(S')} \brkt{ \dbrkt{{\bf G^*}(\phi(e)=S'(e) } -{\bf
d}_{\bf G^*}(S'\ang{e}) } \right| {\bf G^*}(\phi(v))=S'(v),
\forall v\in {\mathbb V}_{1}(S') \right]}^2
\nonumber\\
&=& \brkt{ \Ex_{\phi\in\Phi(h)}\left[\left. \prod_{e\in D} \brkt{
\dbrkt{{\bf G^*}(\phi(e)=S(e) } -{\bf d}_{\bf G^*}(S\ang{e}) }
\right| {\bf G^*}(\phi(v))=S(v),\forall v\in {\mathbb V}_{1}(S)
\right]}^2
\nonumber\\
&=& \brkt{ \Ex_{\phi\in\Phi(h)}\left[\left. \prod_{e\in D} F_e(\phi(e))
\right|
{\bf G^*}(\phi(v))=S(v), \forall v\in
{\mathbb V}_{1}(S)
\right]
}^2\nonumber
\\
& \stackrel{(\ref{070907c})(\ref{lj22b})}\leq & 2\min_{e_0\in D} \Ex_{\varphi\in \Phi(Mh)}
\Ex_{{\bf
e^*}\in {\bf \Omega}_I } [ \brkt{ \Ex_{{\bf e}\in {\bf \Omega}_I
}[F_{e_0}\brkt{\bf e}
 |\,
 {\bf e}\stackrel{
\partial {\bf G^*}/\varphi}{\approx}{\bf e^*}]
}^2 |\,{\bf G^*}(\partial {\bf e^*})=S(\partial e_0) ] \nonumber
\\
&=& 2\min_{e_0\in D} \Ex_{\varphi\in \Phi(Mh)}
\Ex_{{\bf e^*}\in {\bf \Omega}_I } [ \brkt{ \Ex_{{\bf e}\in
{\bf \Omega}_I }[ \dbrkt{{\bf G^*}({\bf e})=S(e_0) } |\,
 {\bf e}\stackrel{
\partial {\bf G^*}/\varphi}{\approx}{\bf e^*}]-{\bf d}_{\bf
G^*}(S\ang{e_0})]
}^2 |\,{\bf G^*}(\partial {\bf e^*})=S(\partial e_0)] \nonumber
\\
&& (\mbox{by the definition of $F_e({\bf e})$ since }\,\,{\bf d}_{\bf
G^*}(S\ang{e})\,\mbox{does not depend on}\,\, {\bf e})
\nonumber\\
& \stackrel{(\ref{061220})}{=}& 2 \cdot \min_{e_0\in D}\eta(S\ang{e_0})
\nonumber\\
& \leq & 2 \cdot \max_{e_0\in \mathbb{V}_2(S)}\eta(S\ang{e_0}).
\label{lj23j}
\end{eqnarray}
Now, choose an $e_0\in {\mathbb V}_2(S)$ which maximizes
$\eta(S\ang{e_0}).$ It then follows from Lemma \ref{lj23c} that
\begin{eqnarray}
&& \Prob_{\phi\in\Phi(h)}[ {\bf G^*}(\phi(e))=S(e),\,\forall e\in
{\mathbb V}_2(S) | {\bf G^*}(\phi(v))=S(v),\,\forall v\in {\mathbb
V}_1(S) ]
\nonumber\\
&\stackrel{L.\ref{lj23c},(\ref{lj23j})}{=}& \prod_{e\in {\mathbb V}_2(S)}{\bf d}_{\bf G^*}(S\ang{e})
\dot{\pm} \sqrt{2}|{\mathbb V}_2(S)| \sqrt{\eta(S\ang{e_0})}
\nonumber\\
&\stackrel{(\ref{071105}),(\ref{lj24m}),(\ref{081021c})}{=}& 
\brkt{ {\bf d}_{\bf G^*}(S\ang{e_0}) \dot{\pm} {\sqrt{2}
{r\choose 2}h^2 \sqrt{\eta(S\ang{e_0})} \over
\brkt{2\sqrt{\epsilon_1}/c_2({\bf G^*})}^{|{\mathbb V}_2(S)|-1} }
} \prod_{e\in {\mathbb V}_2(S),e\not=e_0} {\bf
d}_{\bf G^*}(S\ang{e}) 
\nonumber\\
&\stackrel{(\ref{lj24b})}{=}
& \brkt{ {\bf d}_{\bf G^*}(S\ang{e_0}) \dot{\pm} C
\sqrt{\eta(S\ang{e_0})}}\prod_{e\in {\mathbb V}_2(S),e\not=e_0}
{\bf d}_{\bf G^*}(S\ang{e}) \, \,
(\because c_2({\bf G^*})=c_2({\bf G})\le b_2)
\label{071031a}
\\
&\stackrel{(\ref{lj25})}{=}&\prod_{e\in {\mathbb V}_2(S)} ({\bf d}_{\bf
G^*}(S\ang{e}) \dot{\pm}\delta(S\ang{e})).\quad
\label{lj25a}
\end{eqnarray}
For any $S \in \mathcal{S}_{h,{\bf G^{*}}},$ 
(\ref{lj24}) holds, and we have shown
that $\delta$ satisfies (\ref{lj19}).\\[.1in]
Finally we turn to showing that $\delta$ satisfies (\ref{lj26}).\\
$\bullet$ {\bf [Bounding the average error size] }
For $I\in {\mathfrak{r}\choose
2}$, it follows from the lineality of expectation that 
\begin{eqnarray}
&& \brkt{ \Ex_{n\in [0,\tilde{n}-1],\varphi\in\Phi(m(n))}
\Ex_{{\bf \tilde{e}}\in {\bf \Omega}_I}[
\sqrt{ \eta({\bf G^*}\ang{\bf \tilde{e}})} ] }^2 \nonumber
\\
&\le & \Ex_{n,\varphi} \Ex_{{\bf \tilde{e}}\in {\bf \Omega}_I}[
\eta({\bf G^*}\ang{\bf \tilde{e}}) ]\quad\quad (\mbox{by
Cauchy-Schwarz}) 
\nonumber
\\
&\stackrel{(\ref{061220})}{=}& \Ex_{n,\varphi, {\bf \tilde{e}}}
\Ex_{\varphi^{'}\in\Phi(Mh)}\Ex_{{\bf e^*}\in {\bf \Omega}_I} [
\brkt{\Prob_{{\bf e}\in {\bf \Omega}_I}[ {\bf G^*}({\bf e})={\bf
G^*}({\bf \tilde{e}}) |\, {\bf e}\stackrel{\partial {\bf
G^*}/\varphi^{'}}{\approx}{\bf e^*} ] - {\bf d}_{\bf G^*}({\bf
G^*}\ang{\bf \tilde{e}} )}^2 |\,
 {\bf e^*}
\stackrel{\partial {\bf G^*}}{\approx}{\bf \tilde{e}} ]
 \nonumber\\
&\stackrel{(\ref{081021d})}{=}& 
\mathop{\Ex}_{n,\varphi, {\bf \tilde{e}}} \sum_{\mathfrak{c}_I\in {\rm C
}_I({\bf G^*})}\dbrkt{{\bf G^*(\tilde{e})} = \mathfrak{c}_I}
\mathop{\Ex}_{\varphi^{'}, {\bf e^*}} [ \brkt{\mathop{\Prob}_{{\bf e}\in {\bf
\Omega}_I}[ {\bf G^*}({\bf e})=\mathfrak{c}_I |\, {\bf
e}\stackrel{\partial {\bf G^*}/\varphi^{'}}{\approx}{\bf e^*} ] -
\mathop{\Prob}_{{\bf e}\in {\bf \Omega}_I}[ {\bf G^*}({\bf
e})=\mathfrak{c}_I |\, {\bf e}\stackrel{\partial {\bf
G^*}}{\approx}{\bf \tilde{e}} ] )}^2 |\,
 {\bf e^*}
\stackrel{\partial {\bf G^*}}{\approx}{\bf \tilde{e}} ]
\nonumber \\
& \le & \Ex_{n,\varphi, {\bf \tilde{e}}} \sum_{\mathfrak{c}_I\in
{\rm C}_I({\bf G^*})} \Ex_{\varphi^{'}, {\bf e^*}} [
\brkt{\Prob_{{\bf e}\in {\bf \Omega}_I}[ {\bf G^*}({\bf
e})=\mathfrak{c}_I |\, {\bf e}\stackrel{\partial {\bf
G^*}/\varphi^{'}}{\approx}{\bf e^*} ] - \Prob_{{\bf e}\in {\bf
\Omega}_I}[ {\bf G^*}({\bf e})=\mathfrak{c}_I |\, {\bf
e}\stackrel{\partial {\bf G^*}}{\approx}{\bf \tilde{e}} ] )}^2 |\,
 {\bf e^*}
\stackrel{\partial {\bf G^*}}{\approx}{\bf \tilde{e}} ]
\nonumber
\\
& = & 
\sum_{\mathfrak{c}_I\in{\rm C}_I({\bf G^*})} 
\Ex_{n,\varphi, {\bf \tilde{e}}} 
\left[
\Ex_{\varphi^{'}, {\bf e^*}} 
[
\brkt{\Prob_{{\bf e}\in {\bf \Omega}_I}[ {\bf G^*}({\bf
e})=\mathfrak{c}_I |\, {\bf e}\stackrel{\partial {\bf
G^*}/\varphi^{'}}{\approx}{\bf e^*} ]}^2
 |\, {\bf e^*}\stackrel{\partial {\bf G^*}}{\approx}{\bf \tilde{e}} 
]\right.
 + \brkt{
\Prob_{{\bf e}\in {\bf
\Omega}_I}[ {\bf G^*}({\bf e})=\mathfrak{c}_I |\, {\bf
e}\stackrel{\partial {\bf G^*}}{\approx}{\bf \tilde{e}} ] )}^2
\nonumber
\\
& & \left. 
-
2
\Ex_{\varphi^{'}, {\bf e^*}} 
[
\brkt{\Prob_{{\bf e}\in {\bf \Omega}_I}[ {\bf G^*}({\bf
e})=\mathfrak{c}_I |\, {\bf e}\stackrel{\partial {\bf
G^*}/\varphi^{'}}{\approx}{\bf e^*} ]}
 |\, {\bf e^*}\stackrel{\partial {\bf G^*}}{\approx}{\bf \tilde{e}} ]
 \brkt{
\Prob_{{\bf e}\in {\bf
\Omega}_I}[ {\bf G^*}({\bf e})=\mathfrak{c}_I |\, {\bf
e}\stackrel{\partial {\bf G^*}}{\approx}{\bf \tilde{e}} ] )}
\right]
\nonumber
\\
& = & \sum_{\mathfrak{c}_I\in {\rm C}_I({\bf G})}
\Ex_{n,\varphi, {\bf \tilde{e}}} \left[ 
\Ex_{\varphi^{'}, {\bf e^*}}
 [ \brkt{\Prob_{{\bf e}\in {\bf \Omega}_I}[ {\bf G}({\bf
e})=\mathfrak{c}_I |\, {\bf e}\stackrel{\partial {\bf
G^*}/\varphi^{'}}{\approx}{\bf e^*} ]}^2 |
 {\bf e^*}
\stackrel{\partial {\bf G^*}}{\approx}{\bf \tilde{e}} ] -\brkt{
\Prob_{{\bf e}\in {\bf \Omega}_I}[ {\bf G}({\bf e})=\mathfrak{c}_I
|\, {\bf e}\stackrel{\partial {\bf G^*}}{\approx}{\bf \tilde{e}} ]
)}^2 \right] \nonumber
\\
&&  \mbox{\hspace{2.5in} (since $\stackrel{\partial {\bf
G^*}/\varphi^{'}}{\approx}$ is a refinement of $\stackrel{\partial {\bf
G^*}}{\approx}$) }
 \nonumber\\
&= &|{\rm C}_I({\bf G})|
 \Ex_{\mathfrak{c}_I} \Ex_{n,\varphi, {\bf
\tilde{e}}} \left[ \Ex_{\varphi^{'}\in\Phi(Mh)} [
\brkt{\Prob_{{\bf e}}[ {\bf G}({\bf e})=\mathfrak{c}_I |\, {\bf
e}\stackrel{\partial ({\bf G}/\varphi) /\varphi^{'}}{\approx}{\bf
\tilde{e}} ]}^2 ] -\brkt{ \Prob_{{\bf e}}[ {\bf G}({\bf
e})=\mathfrak{c}_I | {\bf e}\stackrel{\partial ({\bf G}/\varphi)
}{\approx}{\bf \tilde{e}} ] )}^2 \right] \nonumber
\\
&\stackrel{(*)}{\le} &
b_2 \Ex_{0\le n<\tilde{n}} \Ex_{{\bf
\tilde{e}},\mathfrak{c}_I} \left[ \Ex_{\phi^{''}\in
\Phi(m({n+1}))}[ \brkt{ \Prob_{{\bf e}}[ {\bf G}({\bf
e})=\mathfrak{c}_I |\, {\bf e}\stackrel{\partial {\bf
G}/\phi^{''}}{\approx}{\bf \tilde{e}} ] }^2 ] -
\Ex_{{\phi}\in\Phi(m(n))}[ \brkt{ \Prob_{{\bf e}}[ {\bf G}({\bf
e})=\mathfrak{c}_I |\, {\bf e}\stackrel{\partial {\bf
G}/{\phi}}{\approx}{\bf \tilde{e}} ] }^2 ]\right]
\nonumber\\
&= & 
{b_2\over \tilde{n}} \sum_{n=0}^{\tilde{n}-1} \Ex_{{\bf
\tilde{e}},\mathfrak{c}_I} \left[ \Ex_{\phi \in \Phi(m(n+1))}[
\brkt{ \Prob_{{\bf e}}[ {\bf G}({\bf e})=\mathfrak{c}_I |\,
 {\bf
e}\stackrel{\partial {\bf G}/\phi}{\approx}{\bf \tilde{e}} ] }^2 ]
- \Ex_{{\phi}\in\Phi(m(n))}[ \brkt{ \Prob_{{\bf e}}[ {\bf G}({\bf
e})=\mathfrak{c}_I | {\bf e}\stackrel{\partial {\bf
G}/{\phi}}{\approx}{\bf \tilde{e}} ] }^2 ]\right]
\nonumber\\
&= & {b_2\over \tilde{n}} \Ex_{{\bf \tilde{e}},\mathfrak{c}_I}
\left[ \Ex_{\phi\in \Phi(m({\tilde{n}}))}[ \brkt{ \Prob_{{\bf e}}[
{\bf G}({\bf e})=\mathfrak{c}_I |\,
 {\bf e}\stackrel{\partial {\bf
G}/\phi}{\approx}{\bf \tilde{e}} ] }^2 ] -
\Ex_{{\phi}\in\Phi(m(0))}[ \brkt{ \Prob_{{\bf e}}[ {\bf G}({\bf
e})=\mathfrak{c}_I |\, {\bf e}\stackrel{\partial {\bf
G}/{\phi}}{\approx}{\bf \tilde{e}} ] }^2 ]\right]\nonumber
\\
&& 
\hspace{8cm}
\,\, \quad \mbox{(since the sum telescopes!)}
\nonumber\\
&\le&{b_2\over \tilde{n} }
\label{lj20}
\end{eqnarray}
where in (*) above we use Fact \ref{080827a} and 
the property that after $n$ is chosen,
it follows from (\ref{lj24a}) that \\
$m(n+1)\ge Mh+ m(n)
$ and that if $\phi'' ({\mathbb D})\supset 
\varphi({\mathbb D}) \cup \varphi'({\mathbb D})$
(where $\phi''({\mathbb D}), \varphi({\mathbb D}), \phi'({\mathbb D})$
denote the ranges of those functions)
then ${\bf
e}\stackrel{\partial {\bf G}/\phi''}{\approx}{\bf \tilde{e}} $
implies ${\bf e}\stackrel{\partial ({\bf
G}/\varphi)/\varphi'}{\approx}{\bf \tilde{e}}$.
\\[.1in]
\noindent 
Now, for any $I\in {\mathfrak{r}\choose 2},$ it follows from
(\ref{lj25}) that
\begin{eqnarray}
 \Ex_{n<\tilde{n},\varphi\in\Phi(m(n)),{\bf e}\in {\bf \Omega}_I}[\delta({\bf G^{*}\ang{\bf e})}]
&\stackrel{(\ref{lj25})}{\le} & \Ex_{n,\varphi}
\Ex_{{\bf e}\in {\bf \Omega}_I}
[
C \sqrt{\eta({\bf G^{*}\ang{e}})}
+
\Prob_{{\bf e}\in {\bf \Omega}_I}
\left[{\bf G^*\ang{{\bf e}}} \in {\rm BAD}_I \right]\cdot 1
]
\nonumber \\
&= & C\Ex_{n,\varphi} { \Ex_{ {\bf e} }[ \sqrt{\eta({\bf
G^*}\ang{\bf e} )}] } +
\Ex_{n,\varphi}\left[
\Prob_{{\bf e}\in {\bf \Omega}_I}[{\bf G^*}\ang{\bf e}\in
{\rm BAD}_I] \right] \nonumber
\\
& \stackrel{(\ref{lj20})(\ref{lj24m})}{\leq} & { C \sqrt{ b_2\over \tilde{n} } }
+\Ex_{n,\varphi}\left[ \sum_{J\subset I} \Prob_{{\bf e}\in {\bf \Omega}_J}[{\bf
d}_{\bf G^*}({\bf G^*}\ang{\bf e})\le {\sqrt{\epsilon_1}\over
|{\rm C}_J({\bf G^*})|} ] \right] .\label{lj101} 
\end{eqnarray}
However, it is easy to see that for any $\tau > 0$, we
have by the definition of ${\bf d_{G^{*}}}$
\begin{eqnarray}
\Prob_{{\bf e}}\left[{\bf d}_{\bf G^{*}}({\bf G^{*}}\ang{\bf
e})\le \tau \right] = \Prob_{{\bf e}}\left[\Prob_{{\bf e}{'}}
\left[ {\bf G^{*}(e)} = {\bf G^{*}(e{'})} | {\bf G^{*}}(\partial
{\bf e}) = {\bf G^{*}}(\partial {\bf e}{'}) \right]\le \tau \right]
\leq \tau. \label{lj102}
\end{eqnarray}
Hence, using 
(\ref{lj101}) and (\ref{lj102}), we
can write
\begin{eqnarray}
\Ex_{n,\varphi}
\Ex_{{\bf e}\in {\bf \Omega}_I}[\delta({\bf G^{*}\ang{\bf e})}]
&\leq&  { C \sqrt{ b_2\over \tilde{n} } }
+\Ex_{n,\varphi}\left[ \sum_{J\subset I} {\sqrt{\epsilon_1}\over
|{\rm C}_J({\bf G^*})|} ] \right] \nonumber\\
&\stackrel{(\ref{lj24g})}{\leq}& {\epsilon\over 2b_2{r\choose 2}}+
 3
\sqrt{\epsilon_1}  
\nonumber\\
&\stackrel{(\ref{lj24b})}{\leq}& {\epsilon \over b_2 {r \choose 2}}.
 \label{lj24c}
\end{eqnarray}
To show that the expectation of the regularity is small, 
we compute
\begin{eqnarray*}
\Ex_{n,\varphi}[ {\bf reg}({\bf G}/\varphi) ] &\le &
\Ex_{n,\varphi}[ \max_{I\in {\mathfrak{r}\choose [2]}} |{\rm
C}_{I}({\bf G}/\varphi)| \,\Ex_{{\bf e}\in {\bf
\Omega}_I}[\delta({\bf G^*}\ang{\bf e})]]\nonumber
\\
&\stackrel{(\ref{081021b})}{\le} & \Ex_{n,\varphi}[ \sum_{I\in {\mathfrak{r}\choose 2}}
|{\rm C}_I({\bf G}/\varphi)| \,\Ex_{{\bf e}\in {\bf
\Omega}_I}[\delta({\bf G^*}\ang{\bf e})]] \nonumber
\\
&\stackrel{(\ref{lj24c})}{\le} &  \sum_{I\in {\mathfrak{r}\choose 2}} b_2\cdot {\epsilon \over
b_2 {r \choose 2}}\\ &=& \epsilon
\end{eqnarray*}
as required. This completes the proof of Theorem \ref{main}. \qed

\end{document}